\documentstyle[12pt]{article}
\newtheorem{theorem}{\normalsize T\scriptsize HEOREM\normalsize}[section]

\newtheorem{lemma}{\normalsize L\scriptsize EMMA\normalsize}[section]
\newtheorem{corollary}{\normalsize C\scriptsize OROLLARY\normalsize}[section]

\newtheorem{problem}{\normalsize P\scriptsize ROBLEM\normalsize}[section]
\def\endproof{\hfill{$\Box$}} 
\def\dfrac{\displaystyle\frac} 

\title{The Decycling Number of Graphs}
\author{Sheng Bau\thanks{Work supported in part by University of Natal Research Fund URF No. 4509.}\\[3mm]
University of Natal, Pietermaritzburg, South Africa\\[3mm] E-mail: baushe@un.ac.za\\[5mm]}  
\date{}
\begin{document} 
\maketitle 
\begin{center}
\end{center} 

\abstract 
For a graph $G$ and $S\subset V(G)$, if $G - S$ is acyclic, then $S$ is said to be a decycling set of $G$. 
The size of a smallest decycling set of $G$ is called the decycling number of $G$. The purpose of this paper 
is a comprehensive review of recent results and several open problems on this graph parameter. Results to be 
reviewed include recent work on decycling numbers of cubes, grids and snakes(?). 
A structural description of graphs with a fixed decycling number based on connectivity is also presented. 
Graphs with small decycling numbers are characterized. 
\endabstract

\section{The Decycling Number of Graphs}
The minimum number of edges whose removal eliminates all
cycles in a given graph has been known as the {\sl cycle rank\/}
of the graph, and this parameter has a simple expression.
That is, if $G$ is a graph with $p$ vertices, $q$ edges
and $\omega$ components, then the cycle rank (or the Betti number)
of $G$ is $b(G) = q - p +\omega$ (\cite{harary}, Chapter 4).
The corresponding problem for removal of vertices does not
have a simple solution. This latter question is difficult
even for some simply defined graphs.

Let $G = (V, E)$ be a graph. If $S\subseteq V(G)$ and
$G - S$ is acyclic, then $S$ is said to be a {\it decycling set\/}
of $G$. The smallest size of a decycling set of $G$ is
said to be the {\it decycling number\/} of $G$ and is denoted
by $\nabla (G)$. A decycling set of this cardinality is
said to be a {\it minimum decycling set\/}. Determining
the decycling number of a graph is equivalent to finding
the greatest order of an induced forest and the sum of the
two numbers equals the order of the graph.
It was shown in \cite{karp} that determining the decycling number
of an arbitrary graph is $NP$-complete (see Problem 7 on the feedback node set 
in the main theorem of \cite{karp}, which asks for a set $S\subseteq V$ of minimum cardinality 
in a digraph $G$ such that every directed cycle of $G$ contains a member of $S$.). 

Clearly, $\nabla (G) = 0$ if and only if $G$ is a forest, and
$\nabla (G) = 1$ if and only if $G$ has at least
one cycle and a vertex is on all of its cycles. It is also
easy to see that $\nabla (K_p) = p - 2$ and
$\nabla (K_{r, s}) = r - 1$ if $r\leq s$. This is easily
extendable to all complete multipartite graphs. For the Petersen graph $P$, $\nabla (P) = 3$.

We now cite some basic results from \cite{beinekev}.
\begin{lemma} {\rm (\cite{beinekev}, Lemma 1.1)} 
Let $G$ be a connected graph with $p$ vertices, $q$ edges and
degrees $d_1, d_2,\ldots , d_p$ in non-increasing order.
If $\nabla (G) = s$, then
$$\sum_{i = 1}^s (d_i - 1)\geq q - p + 1.$$
\end{lemma} 
\begin{corollary} {\rm (\cite{beinekev}, Corollary 1.2)} 
If $G$ is a connected graph with $p$ vertices, $q$ edges and
maximum degree $\Delta$, then
$$\nabla (G)\geq\frac{q - p + 1}{\Delta - 1}.$$
\end{corollary} 
For graphs regular of degree $r$, one may wonder whether there
is a constant $c$ such that
$$
\nabla (G)\leq\dfrac{q - p + 1}{r - 1} + c?
$$
This is not to be the case, even for cubic graphs (graphs that are
regular of degree $3$). Let $G$ be any cubic graph of order $2n$. Replace each vertex of $G$ with a triangle
and denote the resulting graph by $H$. Then $|H| = 6n$ and
$\nabla (H)\geq 2n$. Thus
$$
\nabla (H) -\dfrac{q - p + 1}{2}\geq 2n -\dfrac{3n + 1}{2}
\geq 2n -\dfrac{3n}{2} = \dfrac{n}{2}.
$$
\begin{problem}
Which cubic graphs $G$ with $|G| = 2n$ satisfy $\nabla (G) = \lceil\frac{n + 1}{2}\rceil$ ?
\end{problem} 
\begin{problem}
Which cubic planar graphs $G$ with $|G| = 2n$ satisfy $\nabla (G) = \lceil\frac{n + 1}{2}\rceil$ ?
\end{problem} 

Let $S\subseteq V(G)$ and define
$$\sigma (S) =\sum_{v\in S} d(v),\quad \epsilon (S) = |E(G|_S)|$$
and denote by $\omega (G)$ the number of components of $G$.
Define the {\it outlay\/} of $S$ to be
$$\theta (S) =\sigma (S) - |S| -\epsilon (S) -\omega (G - S) + 1.$$
\begin{lemma} {\rm (\cite{beinekev}, Lemma 1.3)} 
Let $G$ be a connected graph with $p$ vertices and $q$ edges.
If $S$ is a decycling set of $G$, then $$\theta (S) = q - p + 1.$$
\end{lemma} 
\begin{lemma} {\rm (\cite{beinekev}, Theorem 1.4) } 
If $G$ and $H$ are homeomorphic graphs then $\nabla (G) = \nabla (H)$.
\end{lemma} 
Denote by $\alpha (G)$ and $\beta (G)$ the independence and the
covering numbers of $G$ respectively. Then these two parameters
are related by the equality $\alpha (G) +\beta (G) = |G|$.
\begin{lemma} {\rm (\cite{beinekev}, Theorem 1.5)} 
For any nonnull graph $G$, $\nabla (G)\leq\beta (G) - 1$.
\end{lemma} 
Let $G$ and $H$ be two graphs. Then the {\it cartesian product\/} 
$G\times H$ of $G$ and $H$ is defined by assigning 
$$V(G\times H) = V(G)\times V(H),$$ 
$$E(G\times H) = \{\{(x, y), (x', y')\} : [x = x'\wedge yy'\in E(H)]\vee [y = y'\wedge xx'\in E(G)]\}.$$ 
\begin{theorem} {\rm (\cite{beinekev}, Theorem 1.8)} 
For any graph $G$,
$$2\nabla (G)\leq\nabla (K_2\times G)\leq\nabla (G) +\beta (G).$$
\end{theorem} 
The equalities in Lemma 1.5 are satisfied by some graph of each order.
For example, if $G = K^c_p$, then $\nabla (G) =\nabla (K_2\times G) = 0$
and both equalities hold. Also, for the equality to the lower
bound, if $p\geq 2$ then $\nabla (K_2\times K_p) = 2p - 4 = 2\nabla (K_p)$. 
The path of order $p$ gives equality to the upper bound. 

\section{Cubes}
As we have remarked in the previous section, that the determination of the decycling number of an
arbitrary graph is $NP$-complete \cite{harary}. However, results on the decycling number of
several classes of simply defined graphs have been obtained (\cite{baubdlv}, \cite{baubv} and \cite{beinekev}).

In \cite{beinekev}, upper and lower bounds for the decycling numbers
of cubes and grids were obtained, and in \cite{baubv}, an exact formula
for the decycling numbers of snakes was given.
The results in \cite{baubdlv} and \cite{beinekev} will be reviewed in this section. 

The $n$-dimensional cube (or $n$-cube) $Q_n$ can be defined
recursively: $Q_1 = K_2$ and $Q_n = K_2\times Q_{n-1}$. 
An equivalent formulation, as the graph having the $2^n$
$n$-tuples of $0$'s and $1$'s as vertices with two vertices
adjacent if they differ in exactly one position, gives a
coordinatization of the cube. The following result of \cite{beinekev} 
gives a lower bound on $\nabla (Q_n)$. 
\begin{lemma} {\rm (\cite{beinekev}, Lemma 2.1)} 
Let $n\geq 2$. Then

$(1)$ $\nabla (Q_n)\geq 2\nabla (Q_{n-1})$.

$(2)$ $\nabla (Q_n)\geq 2^{n-1} -\dfrac{2^{n-1} - 1}{n-1}$.
\end{lemma} 
For $n\leq 8$, \cite{beinekev} computed $\nabla (Q_n)$ exactly.
\begin{center} 
\begin{tabular}{c|rrrrrrrr} 
\hline
\hline 
$n$ & 1 & 2 & 3 & 4 & 5 & 6 & 7 & 8\\ 
\hline 
$\nabla (Q_n)$ & 0 & 1 & 3 & 6 & 14 & 28 & 56 & 112\\ 
\hline 
\hline 
\end{tabular} 
\end{center}  
One of the main results of \cite{beinekev} is the computation of
upper and lower bounds of the decycling numbers of
$n$-cubes for $9\leq n\leq 13$. The result is:
\begin{center} 
\begin{tabular}{c|r|r} 
\hline 
\hline 
Cubes & Lower Bound for $\nabla$ & Upper Bound for $\nabla$\\ 
\hline 
$Q_9$ & 224 & 312\\ 
$Q_{10}$ & 448 & 606\\ 
$Q_{11}$ & 896 & 1184\\ 
$Q_{12}$ & 1792 & 2224\\ 
$Q_{13}$ & 3584 & 4680\\ 
\hline 
\hline 
\end{tabular} 
\end{center} 

These results were improved in \cite{baubdlv}. 
\begin{lemma} {\rm (\cite{baubdlv}, Lemma 3.1)} 
For any bipartite graph $G$ with partite sets of cardinality
$r$ and $s$ with $r\leq s$, $\nabla (G)\leq r - 1$.
\end{lemma} 
Since the cartesian product of two bipartite graphs is a
bipartite graph, $Q_n$ is bipartite. With this observation,
the above upper bounds can be lowered a little. For example,
they are $255, 511, 1023, 2047$ and $4095$ respectively.
Applying Lemma 2.1 (2), one can lift the lower bound a little
as well. That is, these lower bounds can be lifted to
$225, 456, 922, 1862$ and $3755$ respectively.
However, one can still go a little further.
\begin{lemma}{\rm (\cite{baubdlv}, Lemma 3.2)} 
If $e$ and $f$ are two adjacent edges of the $n$-cube $Q_n$, then there is a unique $4$-cycle containing
$\{e, f\}$. 
\end{lemma} 
\begin{corollary} {\rm (\cite{baubdlv}, Corollary 3.1)} 

$(1)$ Every edge $uv$ of $Q_n$ is contained in precisely $n - 1$ $4$-cycles;

$(2)$ If $n\geq 3$, then $Q_n$ has precisely $n(n - 1)2^{n-3}$ $4$-cycles. 
\end{corollary} 

Denote by $\rho (u, v)$ the distance between points $u$ and $v$.
Let $x_0\in Q_n$ and define
$$
V_k(Q_n, x_0) = \{x\in Q_n : \rho (x, x_0) = k\}.
$$
Then there is a nice connection between sizes of the sets
$V_k(Q_n, x_0)$ and the binomial coefficients.
\begin{theorem} {\rm (\cite{baubdlv}, Theorem 3.1)} 
$$
|V_k(Q_n, x_0)| = {n\choose k}, \quad 0\leq k\leq n.
$$ 
\end{theorem} 

Let the two partite sets of $Q_n$ be denoted by $X_n$ and $Y_n$.
Then $|X_n| = |Y_n| = 2^{n-1}$. If $x\in X_n$, then the induced
subgraph $G|_{ \{x\}\cup N(x)}$ is a star $S(x)$ of order
$n + 1$ centered at $x$. Call a vertex of degree $1$ of a tree
a {\it leaf\/}. 
\begin{theorem} {\rm (\cite{baubdlv}, Theorem 3.2)} 
For $n\geq 2$, let $x, x'\in X_n$. If $S(x)$ and $S(x')$ are stars
then either $S(x)\cap S(x')=\emptyset$ or $S(x)$ and $S(x')$
have precisely two leaves in common. 
\end{theorem} 
\begin{theorem} {\rm (\cite{baubdlv}, Theorem 4.1)} 

$(1)$ $225\leq\nabla (Q_9)\leq 237;$

$(2)$ $456\leq\nabla (Q_{10})\leq 493;$

$(3)$ $922\leq\nabla (Q_{11})\leq 1005;$

$(4)$ $1862\leq\nabla (Q_{12})\leq 2029;$

$(5)$ $3755\leq\nabla (Q_{13})\leq 4077$.
\end{theorem} 
To obtain the upper bounds given in this theorem, a decyclicng set of the given cardinality is to be 
exhibited in each case and Lemma 2.1 (2) is to be applied. The reader is referred to the elaborate 
proof of this theorem in \cite{baubdlv}. 
\section{Grids}
Another class of graphs for which the decycling number has
been studied to some precision are the grid graphs $P_m\times P_n$.
A standard notation corresponding to matrix notation is
to be adopted for convenience. Thus the $i$th vertex in the
$j$th copy of $P_m$ will be denoted $v_{i,j}$.

If $S$ is a set of vertices in $P_m\times P_n$, then
$S(j)$ will denote the vertices of $S$ in the $j$th column,
and put $S(j, k) = S(j)\cup S(j+1)\cup\cdots\cup S(k)$.
Let $N(j) = |S(j)|$ and $N(j, k) = |S(j, k)|$.

The following results are obtained in \cite{beinekev}.
\begin{theorem} {\rm (\cite{beinekev}, Corollary 3.3)} 
If $m, b\geq 3$, then
$$
\nabla (P_m\times P_n)\geq\left\lfloor\frac{mn - m - n + 2}{3}\right\rfloor.
$$
\end{theorem} 
\begin{theorem} {\rm (\cite{beinekev}, Theorem 5.1)}
For $n\geq 4$,

$(1)$ $\nabla (P_2\times P_n) = \displaystyle\left\lfloor\frac{n}{2}\displaystyle\right\rfloor ;$

$(2)$ $\nabla (P_3\times P_n) = \displaystyle\left\lfloor\frac{3n}{4}\displaystyle\right\rfloor ;$

$(3)$ $\nabla (P_4\times P_n) = n ;$

$(4)$ $\nabla (P_5\times P_n) = \displaystyle\left\lfloor\frac{3n}{2}\displaystyle\right\rfloor -\displaystyle\left\lfloor\frac{n}{8}\displaystyle\right\rfloor - 1;$ 

$(5)$ $\nabla (P_6\times P_n) = \displaystyle\left\lfloor\frac{5n}{3}\displaystyle\right\rfloor ;$ 

$(6)$ $\nabla (P_7\times P_n) = 2n - 1 .$
\end{theorem} 
\begin{theorem} {\rm (\cite{beinekev}, Theorem 5.3)} 
Let $m = 6q + r$ and $n = 6s + t$ with $1\leq r, t\leq 6$.
Then
$$
\nabla (P_m\times P_n)\leq\min\left\{q(2n - 1) +\nabla (P_r\times P_n),
s(2m - 1) +\nabla (P_t\times P_m)\right\}.
$$
\end{theorem} 
\begin{theorem} {\rm (\cite{beinekev}, Theorem 5.4)} 
For $m, n > 2$, 
$$
\nabla (P_m\times P_n) = \frac{mn}{3} + O(m + n).
$$
\end{theorem} 
\begin{theorem} {\rm (\cite{beinekev}, Theorem 5.5)} 
Suppose that $n\equiv 0\quad {\rm (mod 2)}$ and $m = 3r + 1$. 
Then
$$
\nabla (P_m\times P_n) = rn - r + 1.
$$ 
\label{the:gridsd} 
\end{theorem} 
\begin{theorem} {\rm (\cite{baubdlv}, Theorem 5.3)} 
If $S$ is a minimum decycling set of $P_m\times P_n$ with 
\[\nabla (P_m\times P_n) =\left\lceil\frac{mn - m - n + 2}{3}\right\rceil\] 
and 
\[T = \{v_{ij} : i = 2, 4, \cdots , 3m - 2; j = 2, 4, \cdots , 2n - 2\}\] 
then $S\cup T$ is a minimum decycling set of $P_{2m-1}\times P_{2n-1}$. 
\end{theorem} 
\begin{theorem} {\rm (\cite{baubdlv}, Corollary 5.1)} 
For any positive integers $r$ and $s$ 
\[\nabla (P_{6r+1}\times P_{4r-1}) = 8rs - 4r + 1.\] 
\end{theorem} 
This theorem covers some cases other than that covered by Theorem~\ref{the:gridsd}. 

The problem of determining the decycling numbers of the remaining cases of the grid graphs is open.
For the cartesian products $C_m\times C_n$, the following problem is also open.
\begin{problem}
$\nabla (C_m\times C_n) = ?$
\end{problem} 

\section{Snakes}
In this section, a chordless cycle is referred to as a {\it cell\/}. A {\it snake\/} can be defined recursively
as follows. A snake with two cells consists of two cycles with one common edge, one of the two cells will be designated
the {\it head\/} and the other the {\it tail\/}. A snake with $n+1$ cells is obtained from a snake with $n$ cells
by identifying an edge of a new cell with an edge of the tail of the old snake that lies on no other cell. The tail of 
the new snake is the new cell, and the head remains the same. The {\it length\/} of a snake is the number of cells in it.

Determining a minimum decycling set for a snake is algorithmically straightforward. Given a snake $G$,
let $v$ be a vertex on the head with largest possible degree. Put $v$ in the decycling set, then delete it
along with all vertices which lie only on cells that contain $v$. What remains is either a shorter snake or
a single cell. Either repeating this process on the snake which remains (where the new head is the cell
originally adjacent to the old head) or choosing any vertex from a single cell clearly results in a decycling set $S$
of $G$. That $S$ has the minimum possible order follows from the fact that each vertex in $S$ is on some cell
that has none of its other vertices in $S$. Thus, $G$ has a set of $|S|$ vertex disjoint cycles.
Hence $\nabla (G)\geq |S|$.

Let $G$ be a snake. A {\it major pair\/} is a pair of vertices of degree $3$ such that the edge joining them
lies on two cells. A {\it minor pair\/} is a pair of vertices of degree $3$ in a cell which contains
exactly two vertices of degree $4$. A minor pair will be said to lie {\it between\/} the two vertices of
degree $4$. A vertex of degree at least $4$ is called a {\it major\/} vertex.

Note that adding a new tail cell to an existing snake increases the degree of two of its vertices of the old tail
by $1$ each. Since the new tail cell can be incident with at most one vertex of degree at least $3$ in the old snake,
its addition either creates a new major vertex, adds one to the degree of an old major vertex, or creates a major pair.
This idea gives a natural order, from head to tail, $(u_1, u_2, \cdots , u_s)$ to the set of major vertices,
major pairs and minor pairs. Define the {\it name\/} of a snake $G$ to be the sequence $(n_1, n_2,\cdots , n_s)$
where $n_i$ is $3$ if $u_i$ represents a major pair, the degree of the vertex if it represents a major vertex,
and $2$ if it represents a minor pair (where $u_{i-1}$ and $u_{i+1}$ are the vertices of degree $4$ that $u_i$ 
lies between). 

With this definition, there may be several snakes with the same name even if the cells are of uniform length. It is
easily seen that every finite sequence of integers greater than or equal to $3$ is the name of some snake.

Given a snake $G$ and its name $N(G) = (n_1, n_2,\cdots , n_s)$, define the {\it nickname\/} $C(G)$ of $G$ (as a subset of
$\{1, 2,\cdots , s\}$) as follows. 

(1)  $1, s\in C(G)$.

(2) Assume that for $i<s - 1$ it has been determined whether or not each of $1, 2,\cdots i$ is in $C(G)$. Then

(i) If $n_{i+1}\geq 6$, then $i+1\in C(G)$;

(ii) If $n_{i+1} = 5$, then $i+1\in C(G)$ if and only if $i\not\in C(G)$;

(iii) If $n_{i+1} = 4$ then $i+1\in C(G)$ if and only if either $n_i\geq 5$ and $i\not\in C(G)$, or $n_i = 4$
and $i, i-1\not\in C(G)$. 

With this definition, it can be shown that the decycling number of a snake whose cells are all $4$-cycles is the
cardinality of its nickname. From this result the following theorem follows. 
\begin{theorem} {\rm (\cite{baubv}, Theorem 2.1)}
Let $G$ be a snake with nickname $C(G)$. Then $\nabla (G) = |C(G)|$.
\end{theorem} 
A {\it subsnake\/} of a snake $G$ is a subgraph of $G$ that is itself a snake. A {\it straight segment\/} of 
a snake whose cells are all squares is a subsnake in which the vertices of each of the shared edges is a major 
pair. A {\it maximal\/} straight segment $T$ of a square-celled snake $G$ is a straight segment of $G$ such 
that for each cell $s\not\in T$, $T\cup s$ is not a straight segment of $G$. A quare-celled snake $G$ is said to be 
{\it nonsingular\/} if each its maximal straight segment has at least three cells; otherwise it is said to be {\it singular\/}. 
The {\it segment sequence\/} of a square-celled snake $G$ is the sequence of lengths of maximal straight segments 
of $G$ ordered from head to tail. 
\begin{theorem} {\rm (\cite{baubv}, Theorem 3.1)}
If $(d_1, d_2,\cdots , d_k)$ is the segment sequence of a nonsingular snake $G$, then 
\[\nabla (G) =\sum\limits^k_{i=1}\left\lceil\frac{d_i}{2}\right\rceil - k + 1.\] 
\end{theorem} 

The decycling number of a singular snake is certainly related to that of a nonsingular one by means of 
a certain transformation (surgery) of the snake. The decycling numbers of snakes with cell size not equal 
to $4$ are related to those of snakes with cell size $4$ by means of simple transformations. It is therefore 
possible to consider the decycling problem with restriction to square-celled snakes only. 

A snake with triangular cells is a special type of triagulation of a polygon, namely that in which every triangle 
contains at least one edge of the polygon. This raises the question of decycling triangulations of polygons, 
or equivalently, the maximal outerplanar graphs. In general, this seems to be considerably more complicated 
than decycling snakes. At present, we content ourselves with bounds. 
\begin{theorem} {\rm (\cite{baubv}, Theorem 3.4)} 
If $G$ is a maximal outerplanar graph of order $n$, then 
\[1\leq\nabla (G)\leq\left\lfloor\frac{n}{3}\right\rfloor.\] 
\end{theorem} 

Even an algorithm similar to that described at the beginning of this section is not known for the computation of the decycling
number of a triangulation of a polygon. An interesting open problem is to determine the decycling number of outerplanar graphs.
\begin{problem} 
Is there a fast algorithm for computing the decycling numbers of (maximal) outerplanar graphs? 
\end{problem} 
\begin{problem} 
Determine the decycling numbers of $2$-dimensional trees. 
\end{problem} 

\section{Fixed Decycling Numbers}
This section is a simple note on the dependency of decycling number of a graph on its connectivity number. 
In particular, we consider the almost trivial question of determining all graphs with decycling number $2$ or $3$. 

Let $H$ and $J$ be graphs, $S\subseteq V(H)$ and $T\subseteq V(J)$ with $|S| = |T|$. 
Let $f : S\rightarrow T$ be a bijection. An {\it identification\/} of $H$ and $J$ {\it via\/} $f$
is any graph $G$ denoted $H\circ_f J$ obtained by identifying $H$ and $J$ via $f$
and maintaining the adjacencies in the rest of $H$ and $J$. As for the adjacencies in $S$ and $T = f(S)$, 
keep all the edges $E(H|_S)\cup E(J|_T)$. Formally, for
each $x\in S$, identify $x$ and $f(x)$, thus embedding $S$ onto $T = f(S)$. Now for $G = H\circ_f J$,
\[V(H\circ_f J) = V(H)\cup V(J)\backslash T,\] 
\[E(H\circ_f J) = E(H)\cup E(J).\] 
\begin{lemma}
Let $S$ be a $($minimum$)$ decycling set of $H$ and $T$ be a
$($minimum$)$ decycling set of $J$. If $|S| = |T|$ and $f : S\longrightarrow T$ is a bijection, then $S$ is a
$($minimum$)$ decycling set of $H\circ_f J$. 
\end{lemma} 
{\sl Proof\/}: 
Since $S$ and $f(S)$ decycle $H$ and $J$, $S$ is a decycling set of $H\circ_f J$. $H\circ_f J - S$ 
is a forest that is the union of forests $H - S$ and $J - S$. If $S$ is a mimimum decycling set
of $H$ and $T$ is a minimum decycling set of $J$, then $S$ is a minimum decycling set of $H\circ_f J$. 
To see this, let $S'$ be any minumum decycling set of $H\circ_f J$ such that $|S'|<|S|$. Then since 
$\nabla (H) = |S|$, $S'|_H$ cannot decycle $H$. Let $C$ be a cycle of $H - S'|_H$. Then $C$ must be a cycle
of $H\circ_f J - S'$ since $H - S'|_H$ is an induced subgraph of $H\circ_f J - S'$. 
\endproof

Denote by $\kappa (G)$ the connectivity of $G$.
\begin{lemma}
If $\nabla (G) = k$ then $\kappa (G)\leq k + 1$.
\end{lemma} 
{\sl Proof\/}: Let $S$ be a decycling set of $G$ with cardinality
$k =\nabla (G)$. If $\kappa (G)\geq k + 2$, then $G - S$ is
$2$-connected, and hence $G - S$ is not a forest. This contradicts
the choice of $S$. Hence $\kappa (G)\leq k + 1$.
\endproof

Let $G$ be a $(k + 1)$-connected graph with $\nabla (G) = k$.
Since $G$ is not $(k + 2)$-connected, $\kappa (G) = k + 1$.
Let $S$ be a minimum decycling set of $G$. Then $G - S$ is a
connected acyclic graph, i.e., a tree. Thus $G$ can be obtained
by joining a set $S$ to a tree of order at least $k$ in a way
so as to make $G$ a $(k + 1)$-connected graph.
\begin{theorem} 
Let $\nabla (G) = k$. Then $\kappa (G) = k + 1$ if and only if 
there is a tree $T$ and a set $S =\{x : x\not\in T\}$ with $|S| = k$ such that $S$ and $T$ 
are embedded in $G$ and ??. 
\end{theorem} 
{\sl Proof\/}: Let $\kappa (G) = k +1$ and $S$ be a minimum decycling set of $G$. Since $G$ is 
$(k + 1)$-connected, $G - S$ is $k$-connected. Since $S$ is decycling set, $G - S$ is acyclic. 
Hence $G - S$ is a tree. Therefore, $G$ is as described. 

On the other hand, if $G$ is a $(k + 1)$-connected graph as described, then $G - S$ is a tree and 
$\kappa (G) = k + 1$. \endproof 

Inductive construction of all graphs with $\nabla (G) = k$? 

We now determine the graphs with decycling number $2$ and $3$.
\section{Applications} 
The direct application of systems of paths in $n$-cube in coding theory is well known 
(see Chapter ?? of \cite{macwilliamss}), while a maximal such system of paths is what is left after deletion 
of a minimum decycling set. Thus the study of decycling number of cubes, and other particular classes of 
graphs in general, lends itself to coding theory---a more readily applied branch of mathematics. 

\end{document}